\documentclass[12pt]{article} 

\usepackage[latin1]{inputenc}

\usepackage{amscd,amsfonts, amsmath, amsthm, amssymb}
\usepackage{graphicx}
\usepackage{mathrsfs}

\newcommand{\be}{\begin{equation}}
\newcommand{\ee}{\end{equation}}

\newcommand{\cO}{\mathcal{O}}
\title{\bf Singular integration towards a spectrally accurate finite difference operator}
\author{\footnotesize By  ANDR\'E NACHBIN\thanks{Address for correspondence:  Andr\'e Nachbin, IMPA, Estrada Dona Castorina 110, 
Jardim Bot\^anico, Rio de Janeiro, RJ, Brazil, CEP 22460-320; e-mail:  nachbin@impa.br}} 

\newcommand{\fintp}{\int_{0}^{2\pi}\hspace{-2.0em}-  ~  ~   }
\newcommand{\EQ}{\begin{equation}}
\newcommand{\EN}{\end{equation}}

\begin{document}
\date{}
\maketitle

\abstract{

It is an established fact that  a finite difference operator  approximates a derivative with a fixed algebraic rate of convergence. 
Nevertheless, we exhibit a new finite difference 
operator and prove it has spectral accuracy. Its rate of convergence  is not fixed and improves with the function's regularity. 
For example, the rate of convergence is exponential for analytic functions.
Our  new framework is conceptually nonstandard, making no use of polynomial interpolation, nor any other expansion
basis, such as typically considered in  approximation theory.  Our new method arises solely from the numerical manipulation of singular integrals, through an accurate quadrature for 
Cauchy Principal Value convolutions. The  kernel is a distribution which gives rise to multi-resolution grid coefficients.  
The respective distributional  finite difference scheme is spatially structured having stencils of different support widths. 
These multi-resolution stencils test/estimate function variations in a nonlocal fashion, 
giving rise to a highly 
accurate distributional finite difference operator. Computational illustrations are presented, where the accuracy and roundoff error structure 
are compared with the respective 
Fourier  based method. We also compare our method with a recent and popular complex-step method.

%

\section{Introduction and background}
\label{sec_intro}

\par
~

Finite difference operators (FDO) are an important  tool for solving  differential equations in an approximate form. Finite difference methods (FDM) 
are used in computations
as well as in proofs where, for example, a semi-discretization can lead to an existence proof.  Finite difference
approximations typically arise from polynomial interpolation or 
truncating Taylor expansions. In the present work a completely different framework is adopted.

In their classical book \cite{RichtmyerM}, Richtmyer and Morton   
consider an example with the heat equation  $u_t=\sigma u_{xx}$ and provide the convergence rate for the following FDM:
\[
\frac{u_{j}^{n+1} -u_j^n}{\Delta t} = \sigma \frac{ u^n_{j+1}-2u^n_{j}+u^n_{j-1}}{\Delta x^2}, ~~j=1,2,...,J-1; ~~n=0,1,...,
\]
with initial condition 
\[
u_j^0=\varphi(j\Delta x), ~~j=1,2,...,J-1, 
\]
and boundary conditions $u_0^n = u_J^n = 0$, $n=0,1,...$
The interval $[0,L=J\Delta x]$ is partitioned by a uniform grid in space. At a discrete time $t^n=n\Delta t$ the FDM solution at node $x_j=j\Delta x$
 is denoted by $u^n_j$, whereas $u=u(x_j,t^n)$ is the heat equation solution evaluated on the grid.   
For an initial heat profile $\varphi(x)\in C^p$ the following convergence rates, in space and time, are obtained (\cite{RichtmyerM}, page 23):
\begin{equation}
u^n_j - u = \left \{
\begin{array}{l}
\cO(\Delta t^{p/4})  = \cO(\Delta x^{p/2}), ~~\mbox{for} ~~p \leq 3, \\
\cO(\Delta t~|\log \Delta t|)  = \cO(\Delta x^2~|\log \Delta x|), ~~\mbox{for} ~~p=4,\\
\cO(\Delta t)  = \cO(\Delta x^2), ~~\mbox{for} ~~p>4.
\end{array}
\right  .
\end{equation}

Relevant to the present 
work, one observes a convergence rate  partially depending on the regularity of  the initial data.  For a piecewise-linear initial temperature distribution,
the FDM's error decreases  as a square root with respect to the mesh-size. On the other extreme, for an analytic  function
$\varphi(x)$ the error does not decay any faster 
than $\cal{O}$$(\Delta x^2)$.  
The centered difference operator in space, is of second order. 
Irrespective of the importance of FDMs, Richtmyer and Morton  \cite{RichtmyerM} mention that 
{\it  ``In any case, however, finite-difference methods for partial differential equations seldom achieve more than a modest accuracy"}. Fornberg \cite{Fornberg}, among others, 
investigated higher order difference operators for approximating a first derivative.   For a fixed number of grid points $N$, Fornberg  deduced optimal 
coefficients in the operator representation. Fornberg in Table 1 \cite{Fornberg} displays the coefficients for
operators using up to 60 grid points and having an order of accuracy $\cO(\Delta x^p)$, with $p=120$.

In the present work we prove an unexpected result providing the existence of a  finite difference operator with spectral accuracy. 
Namely a discrete operator where the truncation error improves with the regularity
of the underlying function.  As will be shown, our new discrete operator is based on singular integration which generates a family of  second order difference operators, as above. 
Spectral  accuracy is then achieved through (multi-resolution) combinations of these lower order approximations. A striking point is that 
multi-resolution arises from the fact that the kernel of the integral operator is a distribution. 
We do not use a polynomial
approximation as in FDMs, nor  invoke a trigonometric basis as in spectral methods. Our spectrally accurate  FDO arises naturally from the 
numerical manipulation of singular integrals, namely Hilbert transforms.  Harmonic conjugation plays an interesting role as will be discussed
together with the popular (complex-step) method introduced by Squire and Trapp \cite{Squire} and recently advertised by Higham \cite{HighamSN}. The complex-step method aims at reducing the 
truncation error by successfully allowing extremely small step sizes. 
In contrast, in our  method the truncation error will be (conceptually) small in the case of
smooth functions.    Very small step sizes are not needed. The key to this fact is the spectrally accurate quadrature of singular integrals.
\section{Singular integrals and quadrature}

~

Sidi and Israeli  \cite{SidiIsraeli} proved the spectral accuracy of a numerical quadrature for 
periodic singular integrals, defined in the Cauchy Principal Value  sense. 
Spectral accuracy is a remarkable property for a numerical method, where the error depends on the
smoothness of the function. If the function is analytic the error decays exponentially. If the function is band-limited, the 
quadrature is exact up to roundoff error. 
Hence it is also known as exponential convergence or infinite-order accuracy \cite{Canuto}. 
Spectral-accuracy properties  are precisely defined in Fourier space  by combining the Paley-Wiener theorem with
the Poisson summation formula \cite{Trefethen}.

From the work of Sidi and Israeli, theorems 7(a) and 8 play a role in 
our theorem and therefore are respectively reproduced below, in a convenient notation.

\newtheorem{thm}{Theorem}
\begin{thm}  \label{SidiIsraeli7a}
Consider the singular integral $I = \int_{a}^{b}\hspace{-1.3em}- ~~G(x)dx$, where the $L$-periodic integrand has a pole at $x=y$, 
thus being represented in the form $G(x)=g(x)/(x-y)+\tilde{g}(x)$. 
Assume that the functions $g(x)$ and $\tilde{g}(x)$ are $2m$ times differentiable on $[a, b]$. 
Let the grid spacing be $h=L/N$. Letting $y=x_k$, the quadrature given by
\[
Q_N[G]\equiv h \!\!
\sum_{\footnotesize \begin{array}{c} j\! =\! 1\\ j\neq k \end{array} }^{N} G(x_j)
\]
yields the error
\[
E_N[G] = I - Q_N[G] = \left ( g(y)+\tilde{g}(y) \right ) h + {\cal{O}}(h^{2m}),~~\mbox{as}~h\rightarrow 0.
\]
\end{thm}

This quadrature skips the singular point and the error is of first order in $h$. Nevertheless, the next order error-term depends on the regularity of 
the functions $g$ and $\tilde{g}$. Sidi and Israeli  \cite{SidiIsraeli} then used a Richardson extrapolation to take advantage of this fact, 
namely canceling the lower order term and thus obtaining spectral accuracy on the coarser grid. This follows from 

\begin{thm}  \label{SidiIsraeli08}
Consider $G(x)$ and the quadrature $Q_N[G]$  as in Theorem \ref{SidiIsraeli7a}, with $h=h_N\equiv L/N$. Perform the Richardson extrapolation
\[
\tilde{Q}_N[G] = 2Q_{2N}[G]-Q_N[G] = h_N \sum_{j=1}^{N} G(a+jh_N- 0.5 h_N). 
\]
$\tilde{Q}_N[G]$ is a midpoint rule approximation and 
\[
\tilde{E}_N[G] = I - \tilde{Q}_N[G] =  {\cal{O}}(h^{2m}),~~\mbox{as}~h_N\rightarrow 0.
\]
\end{thm}

$\tilde{Q}_N[G]$ is also known as the alternate trapezoidal rule (ATR) as depicted in figure \ref{ATR}. If the integrand's singularity is, for example, at the node
indicated by a circle, then the singular integral is calculated on what we will call the {\it adjoint grid}, indicated by the squares. This leads to spectral accuracy on
a grid of spacing $2h$.
%
%
\begin{figure}[h]
\center
\includegraphics[height=.8in,width=4.5in]{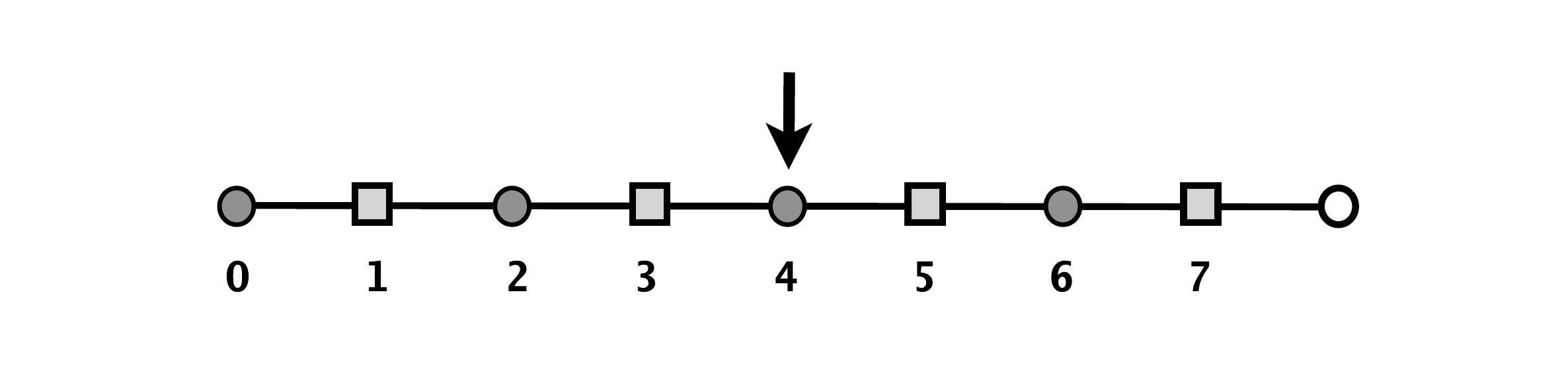} 
\caption{Grid for the alternate trapezoidal rule (ATR). Two alternating grids: one with square nodes and the other with circular nodes.
  For example in the case of a convolution with a distribution,  if the kernel blows up at node 4 as
 indicated by the arrow, only the adjoint grid with square nodes is used.
 The last white node calls attention to an existing periodic condition, namely $0\equiv8$.}
 \label{ATR}
 \end{figure}

 \subsection{The distributional finite difference method}
 ~
 
We define several grid levels labelled by $\ell$, having in mind  multi-resolution. In our study only odd-levels 
$\ell=1,3,\cdots$ will arise. On a resolution level-$\ell$ the grid spacing is $\ell h$
and  the standard second order (centered) finite difference operator  is given by
\EQ
 D^2_{j,\ell} [f] \equiv \frac{f_{j-\ell}-2f_j+f_{j+\ell}}{(\ell h)^ 2}.
 \label{FDM2}
 \EN 
 This discrete operator approximates the second derivative $f^{\prime\prime}(x)$ at the mesh point $x=x_j$. 
 It will become clear why this approximation to a second derivative is relevant in the present study.
 Associated with this discrete operator 
 we have a 3-point stencil, supported over $2\ell+ 1$ nodes indicating its resolution.
 
As in the previous theorems we consider periodic functions.
Without loss of generality we adopt $2\pi$-periodicity. 
The following theorem proves the  spectral accuracy of a conceptually new finite difference operator. 
This result is in contrast with the well 
established property of FDMs having an algebraic rate of convergence. 
Our nonstandard formulation  makes no use of polynomial interpolation, including trigonometric polynomials. 
The new method arises solely from manipulations of singular integrals as presented in the following theorem. 

\begin{thm}  \label{dFDM}																		
Let $f(x)\in\mathcal{C}^m$ be a $2\pi$-periodic function. Consider a grid with nodes spaced by $h=2\pi/N$, where $N$
is a multiple of 4.
In approximating the first derivative of $f(x)$, there exists a finite difference operator with spectral accuracy,  where the error depends on the
function's regularity.  The multi-level  difference  operator $D^1_{\bf ML}[f]_j$ acts on all grid values $[f]$ in approximating the 
derivative  $f^\prime(x_j)$. The operator is defined through 
 \EQ
 f^\prime(x_j) \approx D^1_{\bf ML}[f]_j \equiv
-\frac{2}{N^2}  \sum_{d~ \! \mbox{\footnotesize \bf odd}}^{N/2-1} 
\!\!\ C_d  \sum_{\ell ~ \! \mbox{\footnotesize \bf odd}}^{N/2-1} S_\ell^{-2} \left (
 D_{j+d,\ell}^2  [f]  -
 D_{j-d,\ell}^2 [f] 
 \right ) .
 \EN
The coefficient $C_d$ is related to the distance $dh$ from the target point $x_j$ to the center of the compact stencil defined by the operators
$D_{j\pm d,\ell}^2$. The coefficient $S_\ell$ is related to the resolution level $\ell$. These coefficients have the following expressions:
\be
C_d = \cot ({hd}/{2})
\ee
and 
\EQ
S_\ell = 
  \frac{\sin (h\ell/2)}{h\ell}.
\EN

\end{thm}
\noindent
\textbf{Proof:}                                                                													  
Consider the Hilbert transform on the circle
\EQ
{\cal{H}} [f] (x) \equiv \frac{1}{2\pi} \fintp \cot \left ( \frac{x-y}{2} \right ) f(y)~ dy,
\label{hilb}
\EN
as well as the Hilbert transform of the derivative of $f$, given by
\EQ
{\cal{H}} [f^\prime] (x) \equiv \frac{1}{4\pi} \fintp ~\! \frac{f(x)-f(y)}{\sin^ 2 (\frac{x-y}{2})} ~ dy.
\label{hilbD}
\EN
Using the identity 
$- {\cal{H}}[{\cal{H}} [f] ] (x) = f(x)$, the following expression holds for $2\pi$-periodic functions:
\EQ
f^\prime (x) =  - \frac{1}{8\pi^2} \fintp\!
\cot \left ( \frac{x-y}{2} \right ) \!  \fintp  \left ( \frac{f(y)-f(\tilde{y})}{\sin^ 2 (\frac{y-\tilde{y}}{2} )} \right )
  d\tilde{y}~\! dy.
  \label{hilbRep}
\EN
The ATR, provided by Theorem \ref{SidiIsraeli08},  yields the {\it spectrally accurate} expression
approximating  $f^\prime(x_j)$:
\EQ
f^\prime_j \approx
-\frac{2}{N^2} \!\!\!\! \!\! \sum_{\tiny \begin{array}{c} m\! =\! 0\\(m\!+\!j) \mbox{odd} \end{array}}^{N-1} 
\!\!\!\!\!\!\!\! \cot \left ( \frac{h(j-m)}{2} \right ) 
 \!\!\!\! \!\! \sum_{\tiny \begin{array}{c} n\! =\! 0\\(m\!+\!n) \mbox{odd} \end{array}}^{N-1} 
\!\!\!\!\!\!\! \frac{f_m-f_n}{\sin^ 2 (\frac{h(m-n)}{2})}.
\label{DiffATR}
\EN
The above expression can be recast in the form
 \EQ \label{dFDMa}
 D^1_{\bf ML}[f]_j \equiv
-\frac{2}{N^2}  \sum_{d~ \! \mbox{\footnotesize \bf odd}}^{N/2-1} 
\!\!\ \cot (\frac{h}{2}d)  \sum_{\ell ~ \! \mbox{\footnotesize \bf odd}}^{N/2-1} 
\left ( \frac{(h\ell)^2}{\sin^2 (\frac{h}{2}\ell)} \right )
\left ( D_{j+d,\ell}^2  [f]  -
 D_{j-d,\ell}^2 [f]
 \right ) ,
 \EN
by grouping terms and using the symmetry  and the antisymmetry of the respective kernels,
when $d$ or $\ell$ are greater than $N/2$. In this case these indices are relabelled to the equivalent value
less than $N/2$. 
Due to the periodic grid-structure, care is needed in the above expression:  a renumbering should be performed for the stencil 
indices $P=j\pm d \pm \ell$, when $P$ falls outside  the 
set $\{0,1,2,\cdots,N-1\}$. 
The renumbering is straightforward. 
In the respective difference operator, when $P < 0$ we redefine $P:= N+ P$.
Similarly when $P> N-1$ we redefine $P:= P-N$.  This is in accordance with the numbering
extension indicated by the white squares in the example of figure \ref{MLstencil}.

The distance coefficient  is 
\EQ
C_d = \cot (\frac{hd}{2}),
\EN
and the resolution-level coefficient are expressed as
\EQ
2 S_\ell = 
\left . {\mathscr{S}}_{L}(x) \right | _{x=h\ell} \equiv   \left .  \frac{\sin (\pi x/L)}{(\pi x/L)} \right | _{x=h\ell}; ~L=2\pi.
\EN
This completes the proof. $\Box$

New finite differencing properties are readily available through Theorem \ref{dFDM}. The present differentiation method  
approximating $f^\prime(x)$ is based 
on  numerically integrating  diverse estimates of $f^{\prime\prime}(x)$, obtained through compact stencils of various support widths.   These
{\it multi-scale} estimates are all of (low) $2^{nd}$-order and crude, in the sense that many are evaluated on low-resolution (wide) stencils.  
Nevertheless  it is remarkable that, at the end, the superior spectrally accurate operator is obtained by combining 
multi-resolution stencils centered at all nodes of the adjoint grid.  
These combinations are performed through the resolution-level coefficients $S_\ell$, which are 
grid-point values of a (macroscopic) Whittaker cardinal function  \cite{Whittaker,Trefethen}, which equals 1 at the interval's midpoint
and 0 half a period away.
\begin{figure}
\center
\includegraphics[height=3.2in,width=4.2in]{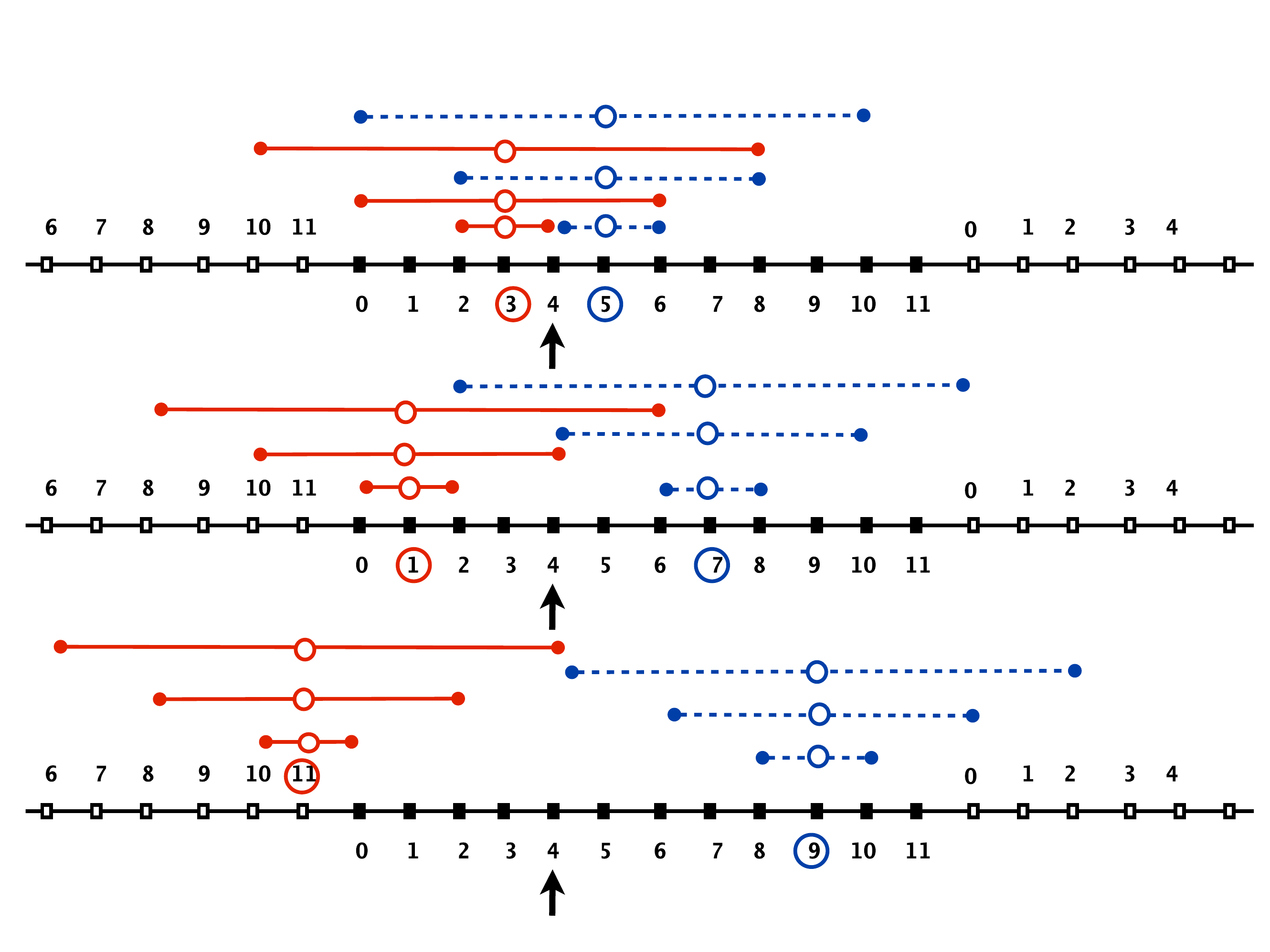}
 \caption{Graphic visualization of the multi-resolution stencils, where $j=4$ and $N=12$.  Dark squares are the mesh points 
 while white squares are their periodic ``extensions", for a better visualization of the multi-level framework.}
 \label{MLstencil}
 \end{figure}

The interpretation from the mathematical, approximation theory,  point of view now follows.  We interpret $f(x)\equiv F(x,y=0)$ as the real-axis trace  of 
a harmonic function $F$ in the upper half-plane.  $F(x,y)$ is the real part of an upper analytic (complex)  function. 
Expression (\ref{DiffATR}) first performs the 
discrete Hilbert transform  of $F_x(x,0)$. 
Despite the fact that differentiation commutes with the Hilbert transform,  expression  (\ref{hilbD})  is exact and
deals directly with increments of $f$ which is very convenient, since values of $f^\prime(x)$ are in principle unknown. 
 Due to the Cauchy-Riemann equations,   this expression can also be viewed as the Hilbert transform of
 $G_y(x,0)$,  where $G$ is  the  harmonic conjugate of $F$. This sets a connection with the interesting work of
 Squire and Trapp \cite{Squire}, as explained below. They  developed a numerical method using complex variables in order to approximate the derivative 
of real valued functions. Their method gained attention \cite{HighamSN} and has been used in many applications, such as numerical Fr\'echet
derivatives  \cite{Higham10} as well as  in Data Science \cite{Deep}. 
As mentioned by Higham \cite{HighamSN} ``{\it 
A fundamental tension in numerical analysis is the interplay between truncation errors and rounding errors; this is particularly prevalent in the numerical approximation of derivatives.}  In order to make finite-difference truncation errors small, one needs to compute with very small values of $h$. This choice might influence
the buildup of  roundoff error and the complex-step method \cite{Squire} is well suited to avoid this buildup. 
In our case truncation-errors might  already be small due to spectral accuracy. Hence our discretization parameter $h$
does not need to be very small, as will be numerically illustrated. 

The complex-step method of Squire and Trapp \cite{Squire} 
extends $f$ as a complex analytic function, by virtue of $f(x_j+ih)$, which by Taylor expansion gives 
\[
F(z)\equiv f(x_j+ih) = f(x_j) + ih f^\prime (x_j) - 0.5 h^2 f^{\prime\prime}(x_j)  + \mbox{$\cal{O}$}(h^3).
\]  
Their second order approximation arises as
\EQ
f^\prime (x_j) \approx \mbox{Im}(f(x_j+ih))/h.
\label{ComplexStep}
\EN
This complex extension of $f$ automatically implies that the respective
harmonic conjugate is zero at $x_j$. This fact avoids a  subtractive cancellation error \cite{Squire}  present in 
FDMs when $h$ is extremely small.
The complex-step expression (\ref{ComplexStep}) is actually the  finite 
difference  approximation of the $y$-derivative of the harmonic conjugate, about the real axis. Namely 
 \EQ
f^\prime (x_j) \approx \frac{\mbox{Im}(f(x_j+ih))-\mbox{Im}(f(x_j))}{h}
\EN
is a semi-discrete approximation of a Cauchy-Riemann equation along the real axis.
In principle this method needs a complex evaluation tool, where care is needed \cite{Higham10}. Our spectrally accurate method 
presented in Theorem \ref{dFDM} 
also made use of the  $y$-derivative of the harmonic conjugate, but without any use of complex functions nor any  Taylor expansion. 
Note that the continuation into the complex plane is different between our method and the complex-step method.
Each method generates a different  analytic continuation in the upper-half plane where  real-axis traces, of  the real part, both coincide with $f(x)$.

 It is instructive to consider a small grid example in order to write explicitly all terms  of our numerical scheme. This allows a clear visualization of all 
 stencils activated through the (distribution) kernels. Consider the example with $N=12$ and take as target node $j=4$.
 Using expression  (\ref{DiffATR}) or (\ref{dFDMa}),  we have that:
 \[
 f_4^\prime \approx   -\frac{2}{N^2}  \cdot
 \left \{  
C_1 \left ( \frac{D_{5,1}^2 - D^2_{3,1}}{S_1} +    
 \frac{D_{5,3}^2 - D^2_{3,3}}{S_3} +
 \frac{D_{5,5}^2  - D_{3,5}^2}  {S_5} 
  \right . \right )  + 
 \]  
 \[
 + C_3
\left (    \frac{D_{7,1}^2  - D_{1,1}^2}{S_1} +   \frac{D_{7,3}^2- D_{1,3}^2}{S_3} 
\!\! +
  \frac{D_{7,5}^2 - D_{1,5}^2} {S_5} \right ) \! + \!
 \]  
 \[
 \left .
+ C_5 \left ( \frac{D_{11,1}^2  - D_{9,1}^2}{S_1}  +
  \frac{D_{11,3}^2 - D_{9,3}^2}{S_3}    +
   \frac{D_{11,5}^2 - D_{9,5}^2}{S_5} \right  ) \!\!
  \right \} \!\!.
 \]  
The stencil notation $D^2_{k,\ell}$, as in(\ref{dFDMa}), has the first subscript $k$ indicating the stencil's center-node, while the second 
subscript equals half of the stencil's 
width.
The multi-resolution stencils are  depicted in figure \ref{MLstencil}. 

These stencils, of varying support width,  explore the entire domain  always centered on the 
adjoint grid. 
This structure is a consequence of convolutions with 
distributions. We will refer to our method as a distributional finite difference method (dFDM),
where stencils of variable compact support test/probe the function's  derivative on the entire domain, resulting in approximations of superior accuracy. 

 \subsection{The differentiation matrix}
 ~
 
 The present method and, for example, the Fourier spectral method have the same spectral accuracy. A natural question is whether
 these methods are identical, where ours operates in the physical domain while the Fast Fourier Transform (FFT) operates in the frequency domain.
 This is answered by looking at the differentiation matrix of each formulation and by examining,  in the next section,  the respective roundoff error 
 structure. 
 
 The standard centered finite difference is obtained from 
 \be  \label{Taylor}
 f^\prime(x_j) = \frac{f(x_{j+1})-f(x_{j-1})}{2h} +\frac{h^2}{12} f^{\prime\prime\prime}(\xi), ~~~x_{j-1}<\xi<x_{j+1} .
 \ee
 The differentiation matrix is trivially obtained as
 \be  \label{diffFDM}
\mathcal{D}{^N}_{\!\!\!\\\!\!\mbox{\bf \tiny FDM}} =   \frac{1}{h}~
\left [
\begin{array}{ccccccc}
0 & 1/2 &   & &  \cdots& &  -1/2\\
  -1/2 & 0 &  1/2 & &   & & \\
  &  -1/2&  0 &   1/2 &  &  &\\
&  & & 0 &   & &   \\
\vdots& & & & \ddots & & \vdots\\
\vdots &  &  &  &    & 0 & 1/2\\
1/2 &\cdots  & & & &-1/2 & 0
\end{array}
\right ] .
\ee

 The differentiation
 matrix for the Fourier method is obtained from interpolating $f(x)$ with a basis of band-limited functions \cite{Trefethen}.
 On the real line  $f(x)$ is approximated by a linear combination of 
translated Whittaker cardinal functions \cite{Whittaker} $\mathscr{S}_h(x)$:
\EQ
p(x) = \sum_{m=-\infty}^{\infty} f_m {\mathscr{S}}_h(x-x_m),
\EN
where ${\mathscr{S}}_h(x)=\sin(\pi x/h)/(\pi x/h)$ are the band-limited functions, whose Fourier spectrum runs exactly up to the Nyquist frequency.
The derivative of $f$ is approximated as
\EQ
f^{\prime}_j \approx \frac{dp}{dx}(x_j) = \sum_{m=-\infty}^{\infty} f_m {\mathscr{S}}^\prime_h(x-x_m).
\EN
In the $2\pi$-periodic case a similar interpolation is made, but with periodic sinc functions \cite{Trefethen}, where
\EQ
p(x) = \sum_{m=1}^{N} f_m \tilde{\mathscr{S}}_N(x-x_m), ~~\mbox{with}~~ \tilde{\mathscr{S}}_N(x)\equiv \frac{\sin(\pi x/h)}{(2\pi x/h)\tan(x/2)},
~N=2\pi/h.
\EN
The functions $\tilde{\mathscr{S}}_N(x)$ are  the periodic counterpart of Whittaker's cardinal functions.
The derivative,  evaluated at grid-points,  equals
\EQ
\tilde{\mathscr{S}}^\prime_N(x_j) = \left \{ 
\begin{array}{ll}
0, & j = 0, \\ 
\frac{1}{2} (-1)^j \cot(jh/2), & j\neq 0.   
\end{array}
\right .
\EN
This leads to the differentiation matrix 
%
%
\be  \label{diffM}
\mathcal{D}^N_{\mbox{\bf \tiny FFT}} =   \frac{1}{2}~
\left [
\begin{array}{ccccccc}
0 & \cdots &   & &  \cdots& &  - C_1\\
  - C_1 & 0 &   & &   & & + C_2\\
 + C_2 &  &  0 &   &  &  & -C_3\\
-C_3&  & & 0 &   & &  +C_4 \\
\vdots& & & & & & \vdots\\
\vdots &  &  &  &    & \ddots & C_1\\
C_1 &\cdots  & & & \cdots& & 0
\end{array}
\right ] ,
\ee
where we recall that
\(
C_j = \cot ( {jh}/{2} ).
\)
 
The dFDM differentiation matrix is found by recalling that 
the ATR  uses two grids.  
To obtain the distributional finite-differencing matrix it is convenient to rewrite expression (\ref{DiffATR})   into two components, one for 
each grid. Relabelling accordingly yields
\be
 \label{dFD_grid1}
 f^\prime(x_j) \approx 
-\frac{2}{N^2}  \left \{  \alpha_N \!
\sum_{m_1}^{N-1} 
\!\!  \left [ \cot (\frac{h}{2}(j-m))  \right ] ~f_{m_1}
\right .  -
\ee
\be
 \label{dFD_grid2}
\left . 
- \sum_{m_2}^{N-1} 
\left [  \sum_{\ell = 1,~ \! \mbox{\footnotesize \bf odd}}^{N/2-1} 
 \frac{ \cot (\frac{h}{2}(j-m-\ell))+ \cot (\frac{h}{2}(j-m+\ell))} {\sin^2 (\frac{h}{2}\ell)} \right ] f_{m_2} \right \} ,
\ee
where
\be
\alpha_{N}\equiv \sum_{\ell=1, ~ \! \mbox{\footnotesize \bf odd}}^{N/2-1} 
\frac{2}{\sin^2 (\frac{h}{2}\ell)} .
\label{alphaN}
\ee
If the index  $j$ of the target node is even, then indices $m_1$ are odd-valued while indices $m_2$ are even. If $j$ is odd
then  $m_1$ runs over the even grid while indices $m_2$ run over the odd grid. 
The number 2, which appears  in the multi-level coefficient $\alpha_{N}$, is reminiscent of contributions from the
center-nodes of the respective stencils. As mentioned earlier,  the center of the test-stencils are on the adjoint grid. The diagonal entry,
corresponding to $m_2=j$, is zero. The same is true for the entry ($j,m_2$), $m_2=j\pm N/2$, half a period away. 
%

  
\begin{figure}[h]
\center
\includegraphics[height=2.in,width=5.in]{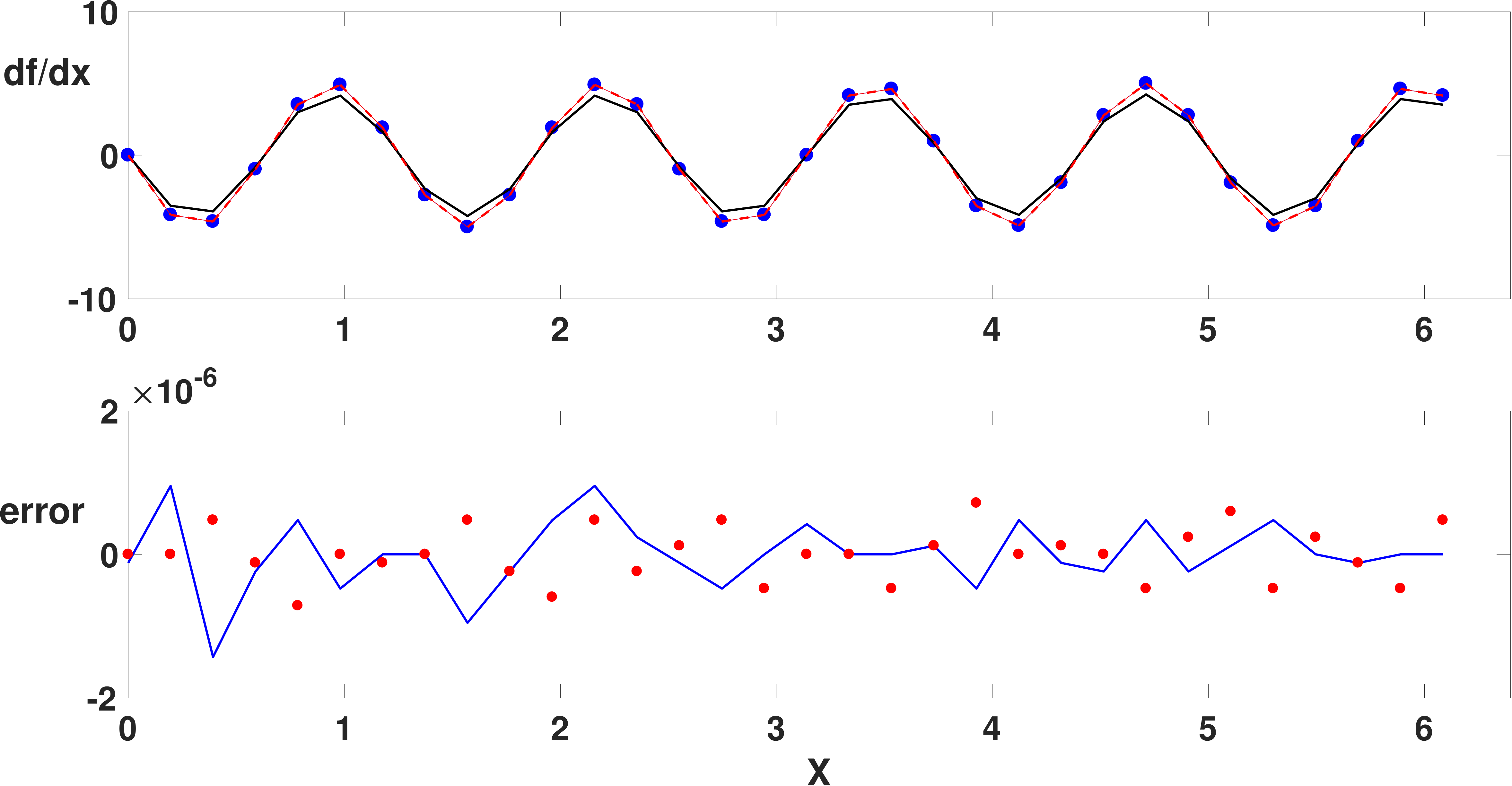}
 \caption{Numerical illustration for $f(x)=\cos(5x)$, with $N=32$.  Top:   Dots display the result from the dFDM. 
Two lines are superimposed connecting the dots: 
 results from the FFT method and the exact derivative. The solid line with a small discrepancy displays the result from the FDM. Bottom: roundoff error 
 difference (between exact $f^\prime$ and method) for
 the dFDM (solid line) and the FFT (dots).
The $\ell_\infty$ norms are: $||f^\prime\mbox{-dFDM}||_\infty=1.4\cdot 10^{-6}$ and $||f^\prime\mbox{-FFT}||_\infty=7.1\cdot 10^{-7}$.}
 \label{cos5}
 \end{figure}

\begin{figure}[h]
\center
\includegraphics[height=2.in,width=5.in]{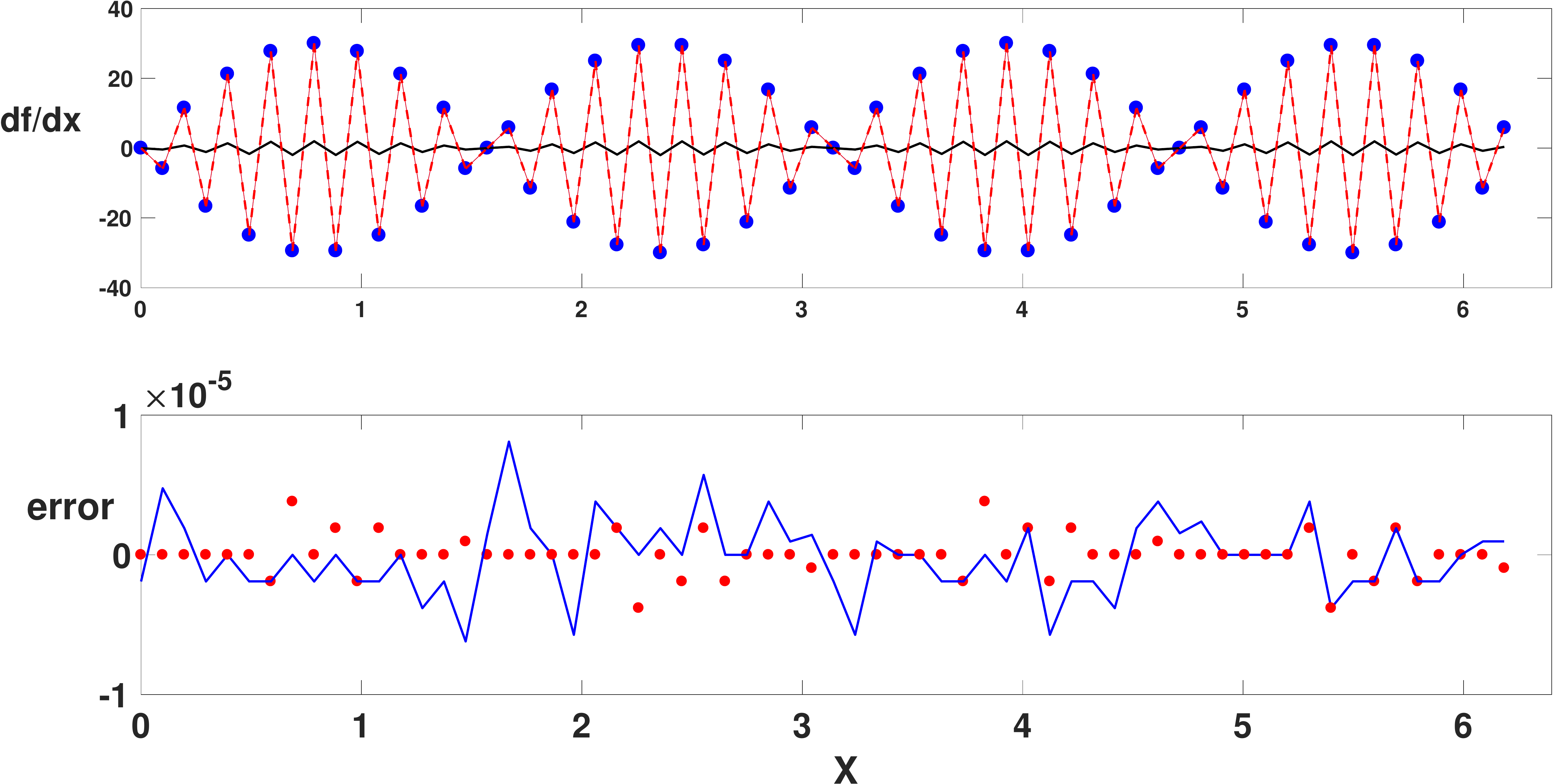}
 \caption{Numerical illustration for $f(x)=\cos(30x)$, with $N=64$. The rapidly varying function is almost at the Nyquist frequency.
 Top: dots display the result from the dFDM. Two lines are superimposed connecting the dots: 
 results from the FFT method and the exact derivative. The solid line with a large discrepancy displays the result from the FDM. This is due to large values
 of the third derivative in (\ref{Taylor}).
 Bottom: the roundoff error difference (between exact $f^\prime$ and method) for
 the dFDM (solid line) and the FFT (dots). 
 The $\ell_\infty$ norms are: $||f^\prime\mbox{-dFDM}||_\infty=8.1\cdot 10^{-6}$ and $||f^\prime\mbox{-FFT}||_\infty=3.8\cdot 10^{-6}$.}
 \label{cos30}
 \end{figure}


\section{Computational illustration}
   ~
   
As our first computational illustration we consider $f(x)=\cos(5x)$ with only $N=32$ points. Due to spectral accuracy, the dFDM and the FFT method are
exact up to roundoff error. The FDM is of second error. At the top of  figure \ref{cos5} we have the derivative computed by these three methods compared to the
exact derivative. The dots are values from the dFDM, coinciding with high accuracy with the FFT and the exact solution. In order to better visualize the roundoff
error all computations were done in MATLAB using single precision. The roundoff error difference between the FFT and the exact derivative are depicted by dots, in  the bottom
part of figure \ref{cos5}, while the roundoff error difference
between the dFDM and the exact value by a solid line. In the top graph we observe the result with the FDM, in a dark solid line, exhibiting some discrepancy with respect
to the other 
methods.

Next we consider a rapidly varying function: $f(x)=\cos(30x)$. We use a grid with $N=64$ points.  We have a well resolved Fourier mode which is close
to the  limiting  Nyquist frequency. Again the values from our dFDM, the FFT method and the exact derivative coincide with high accuracy. 
At the top of figure \ref{cos30} we display these solutions while at the bottom graph the respective roundoff error. At the top graph of figure \ref{cos30} we see a low
amplitude jagged curve representing the FDM approximation which is very poor. This is expected because, as indicated by (\ref{Taylor}), the FDM's 
truncation error scales like $\varepsilon^{-3}$ for rapidly varying functions $f(x/\varepsilon)$, $\varepsilon\ll 1$. 

\begin{figure}[h]
\center
\includegraphics[height=2.5in,width=5.in]{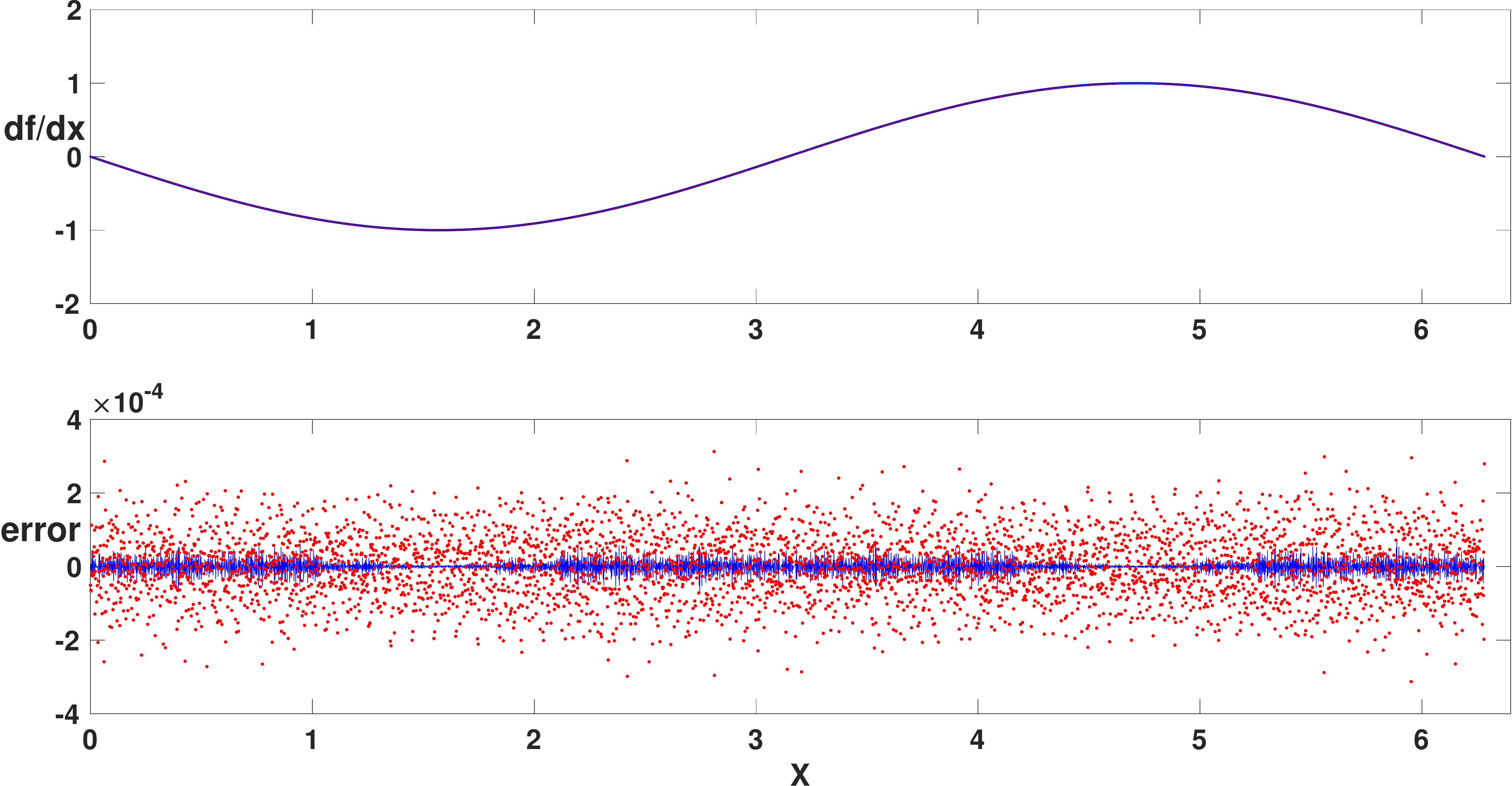}
 \caption{Numerical illustration for $f(x)=\cos(x)$, with $N=4096$. 
Top: Four lines are superimposed: 
 results from the dFDM, the FFT method, the FDM and the exact derivative. We have a slowly varying function and test the methods in the 
 presence of a large number of grid points. The agreement is very good.  Bottom: the roundoff error difference (between exact $f^\prime$ and method) for
 the dFDM (solid line) and the FFT (dots). 
 The $\ell_\infty$ norms are: $||f^\prime\mbox{-dFDM}||_\infty=7\cdot 10^{-5}$ and $||f^\prime\mbox{-FFT}||_\infty=3.1\cdot 10^{-4}$.
  In figure \ref{cosDet} we take a closer look.}
 \label{cos}
 \end{figure}

\begin{figure}[h]
\center
\includegraphics[height=1.5in,width=4.0in]{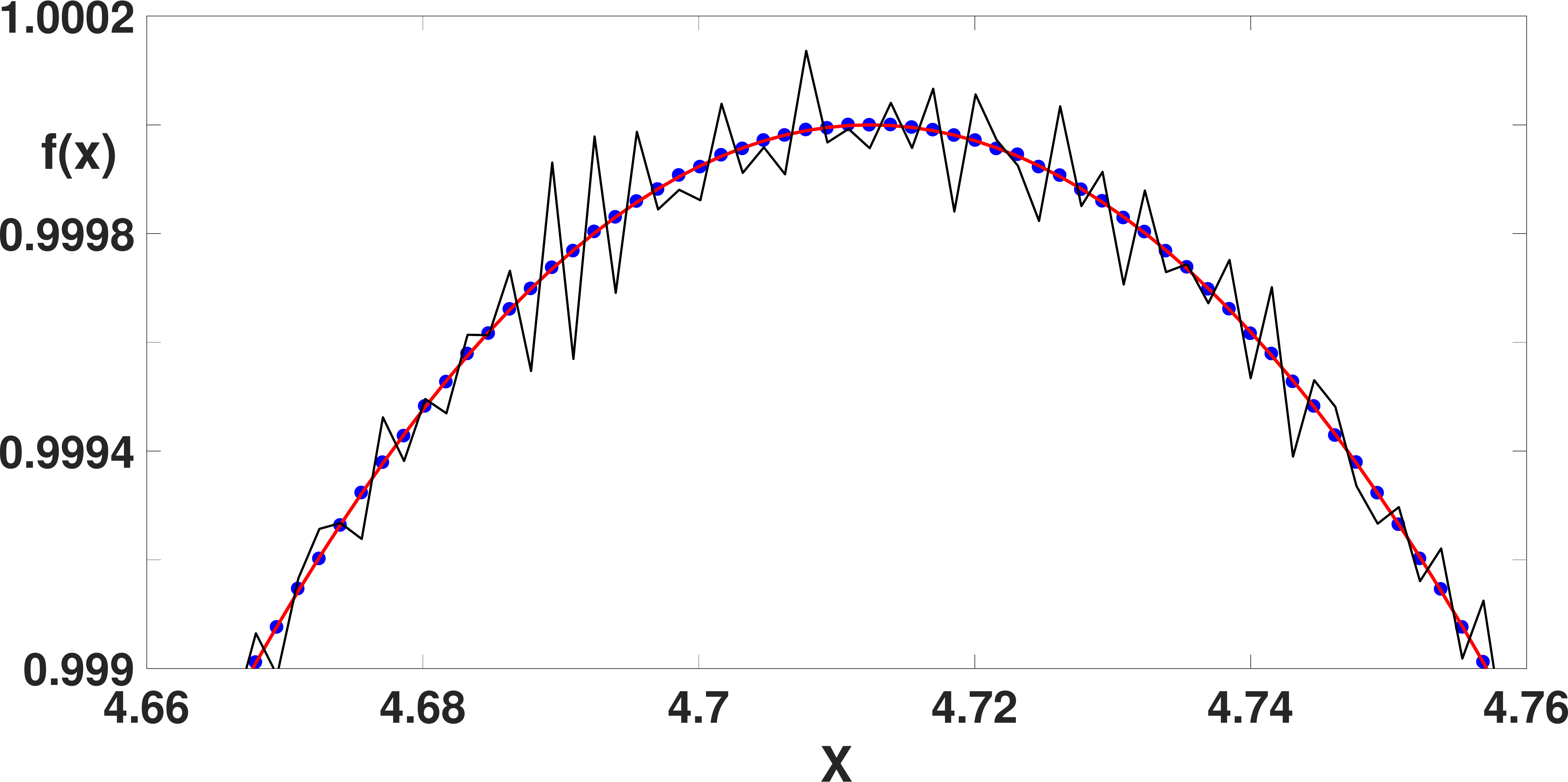}
 \caption{Detail from figure \ref{cos}. Dots display the result form the dFDM which coincide with the exact derivative.
 The sawtooth profile regards the FFT method, where Fourier amplitudes are multiplied by $ik$. The Nyquist frequency,
 containing only roundoff error is amplified by $N/2=2048$ and is visible as we zoom into the graph. Nevertheless this is not visible on the 
 original scale depicted in figure \ref{cos}.}
 \label{cosDet}
  \end{figure}
  
In some problems it might be necessary to work with a great number of grid points.  In figure \ref{cos} we have the derivative of 
a slowly varying function given by $f(x)=\cos(x)$. In this case the truncation error is zero but we test the use of a large number of
grid points:  $N=4096$. The top graph
 in figure \ref{cos} has three curves, regarding the exact derivative, the result of the dFDM and of the FFT method. The roundoff error difference
 with the exact solution is displayed at the bottom. The cloud of points regards the FFT method while the solid line regards the dFDM.  The roundoff error
 structure is quite different, where for the dFDM it is non-uniformly distributed. To have a closer look, refer to
 figure \ref{cosDet} where we have zoomed into the neighborhood of $x=3\pi/4$. We observe that the  FFT method 
 displays oscillations due to the amplification of the sawtooth mode, namely at the Nyquist frequency. 
 Differentiation in the FFT  method
 is performed through its Fourier multiplier, where Fourier amplitudes are multiplied by $ik$.  In the case of a smooth function the upper 
 (inactive) end of the Fourier spectrum is amplified by an order of $10^3$. The agreement between the dFDM and the exact derivative is excellent.

 Now we consider a function with a full Fourier spectrum.  In Fourier space the Gaussian $f(x)=\exp(-(x-\tilde{x})^2/\delta)$ is also a Gaussian.
 We choose  $\delta$ so that $f$ decays fast enough to numerically be well approximated by its periodic extension.  
 At the top of figure \ref{Gauss512} we display  
 the function $f(x)$ while the middle graph 
 superimposes three curves:  the exact derivative, and the results for the dFDM and the FFT method. A total of 
 $N=512$ points were used. The roundoff error is displayed at the bottom part of figure \ref{Gauss512}. The cloud of points regards the 
 roundoff error difference between 
 the exact derivative and the FFT method. The solid line regards the roundoff error difference from the dFDM. The structure is quite different. While the FFT produces 
 roundoff errors uniformly along grid points, the dFDM produces roundoff error only on points where the function $f$ is effectively supported. In figure \ref{Gauss2048} 
 we increase the resolution by using $N=2048$ and the structure is unchanged.
  
The next illustration considers the non-periodic function used by Squire and Trapp \cite{Squire},  in applying their complex-step method.  
They computed the derivative of $f(x)=x^{9/2}$ at the point $\tilde{x}=1.5$. Double-precision was considered, so we now adopt double-precision in Matlab. 
Since our method applies to periodic functions we used a smooth cutoff in the form of a super-Gaussian $sG(x)$ given by
\be
sG(x) = e^{-(\sigma (x-\pi))^s}.
\ee
 To compute the dFDM exactly at $\tilde{x}$  we perform the shift $f(x)=(x-\pi+1.5)^{9/2}$ and therefore
 centered the super-Gaussian at $x=\pi$. Since our method is non-local we explored with different widths of this ``table-top"
 function. We avoid having its effective support in the region where we have the square root of a negative function. 
 We apply the dFDM to the function $\tilde{f}(x)\equiv f(x)\cdot sG(x)$ in order to compute the derivative at $x=\pi$, the shifted position 
 with respect to $\tilde{x}$ of 
 Squire and Trapp.  Excellent results are obtained. In figure \ref{sGauss10_512} we
 used a super-Gaussian with $s=10$, $\sigma = 1.6$ and $N=512$, which yields $h\approx0.012$. 
 The top graph displays the super-Gaussian centered at $x=\pi$. The middle graph displays (in dots)
 the dFDM solution while in  gray we have the exact derivative. The bottom graph depicts the function $f(x)$.  
 At the point of interest the exact derivative value, evaluated through Matlab, is $f^\prime=18.600812734259758$. The approximate value 
 computed with the dFDM is
 $f^ \prime_{257}\approx 18.600812734259215$.  The dFDM value provides 14 digits of accuracy. 
 
 Squire and Trapp \cite{Squire} compare, in Table 1 (page 111), the accuracy of their complex-step method  (CSM)
 with the central finite differencing (FD) $f^\prime\approx f(1.5+h)-f(1.5-h)/(2h)$. 
 Both methods are second order accurate. We reproduce some entries from their table. For a resolution of $h=0.01$, which corresponds to their first table entry, 
 the FD gives $f^\prime\approx 18.602018344501897$ while for the CSM they obtain
 $f^\prime\approx 18.599607128036329$. Rounding at the first 4 digits agrees with the exact value. To obtain 14 digits of accuracy
 the CSM  requires   $h=10^{-7}$ and gives $f^\prime\approx 18.600812734259637$. The dFDM obtained 14 digits with $h=10^{-2}$.
 
 We further reduced the resolution to $N=128$  ($h\approx 0.05$) using $s=8$ and $\sigma =3.0$. These choices  provide a sharp super-Gaussian 
 cutoff as shown by the solid line in figure \ref{sGauss4_8_128}.
 We obtain $f^\prime_{65}\approx   18.613734049360190$, namely with 3 digits of accuracy.
 Due to the low resolution we then explored with a smoother super-Gaussian, adopting $s=4$ and $\sigma=1.6$, as depicted by the dotted line in figure 
  \ref{sGauss4_8_128}.
 We obtain $f^\prime_{65}\approx  18.600812731981232$, namely with 10 digits of accuracy. The 
 different super-Gaussians were tested in an
 experimental fashion. We have not  attempted to find it's optimal width.
  
   %
\begin{figure}[t]
\center
\includegraphics[height=2.5in,width=5.5in]{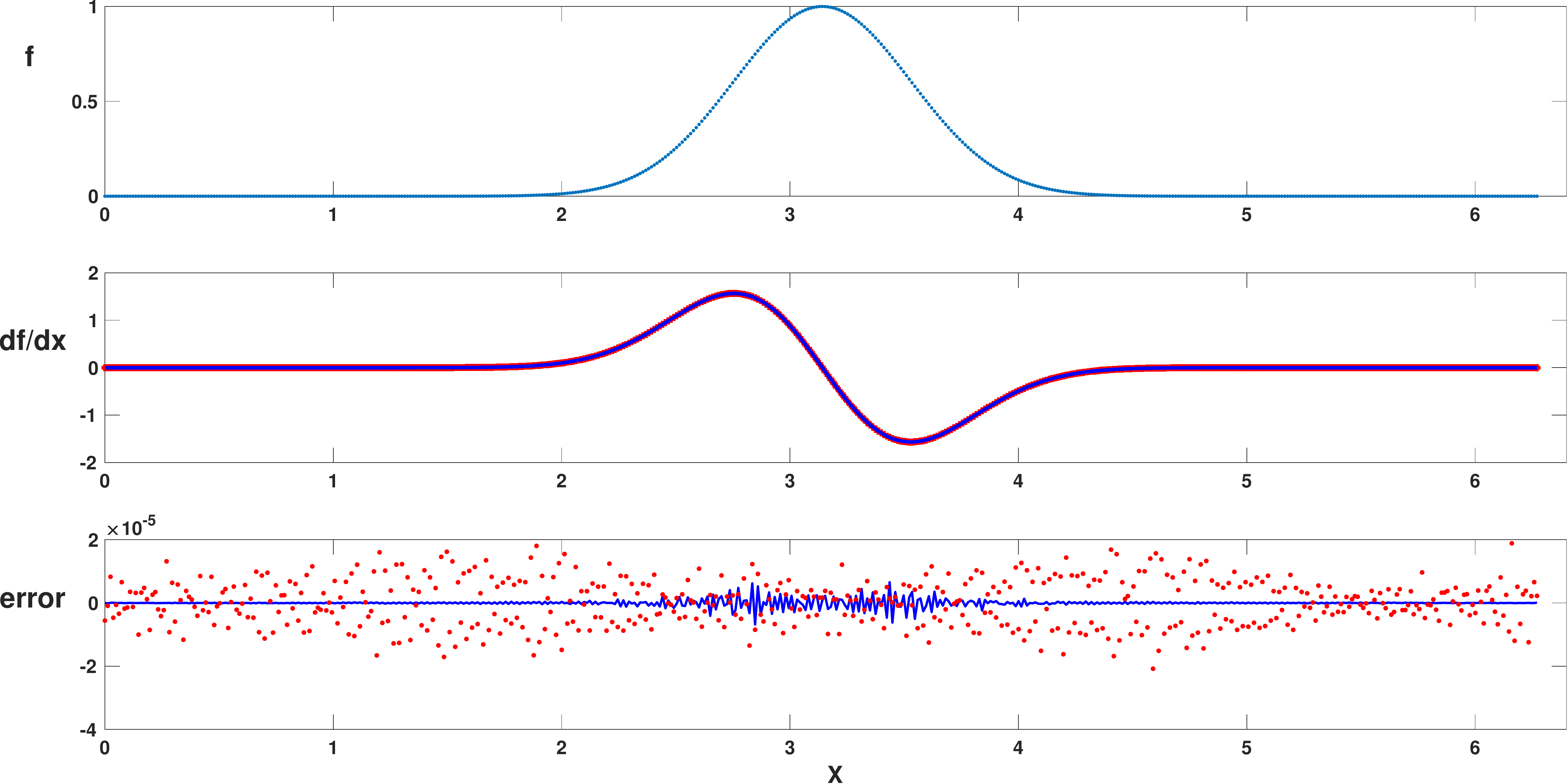}
 \caption{Top: the Gaussian $f(x)=\exp(-(x-\pi)^2/0.3)$. The grid has $N=512$ points. Middle: three curves coincide, 
 namely the exact derivative, the dFDM and the FFT method. Bottom: the roundoff error difference for the dFDM (solid line) and 
 the FFT method (dots). The $\ell_\infty$ norms are: $||f^\prime\mbox{-dFDM}||_\infty=7.8\cdot 10^{-6}$ and $||f^\prime\mbox{-FFT}||_\infty=2.1\cdot 10^{-5}$.}
 \label{Gauss512}
 \end{figure}
\begin{figure}[t]
\center
\includegraphics[height=3.0in,width=5.5in]{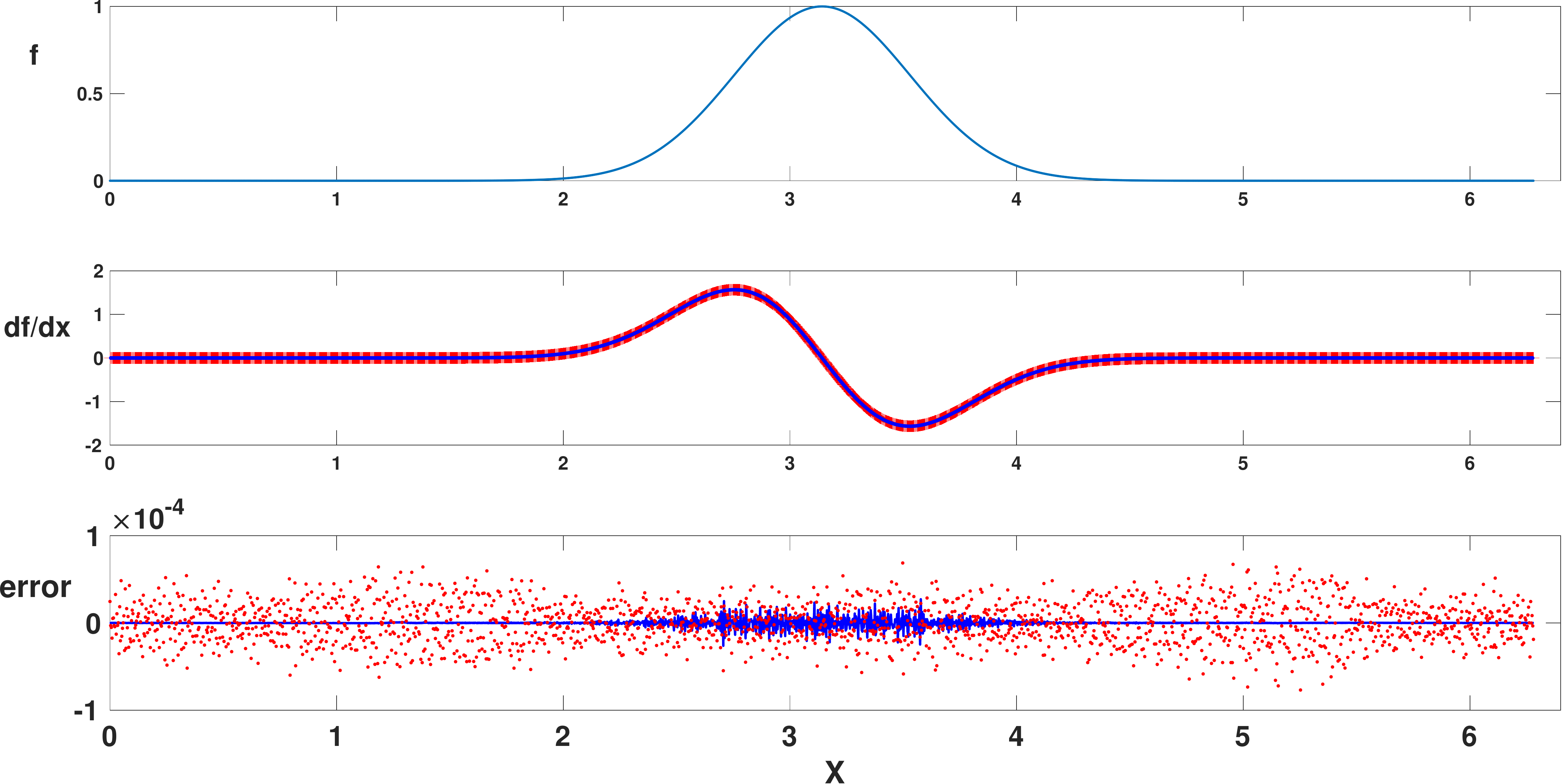}
 \caption{Top: the Gaussian $f(x)=\exp(-(x-\pi)^2/0.3)$. The grid has $N=2048$ points. Middle: three curves coincide, 
 namely the exact derivative, the dFDM and the FFT method. Bottom: the roundoff error difference for the dFDM (solid line) and 
 the FFT method (dots). The $\ell_\infty$ norms are: $||f^\prime\mbox{-dFDM}||_\infty=2.7\cdot 10^{-5}$ and $||f^\prime\mbox{-FFT}||_\infty=7.7\cdot 10^{-5}$.}
 \label{Gauss2048}
 \end{figure}
 
\begin{figure}[t]
\center
\includegraphics[height=2.5in,width=5.5in]{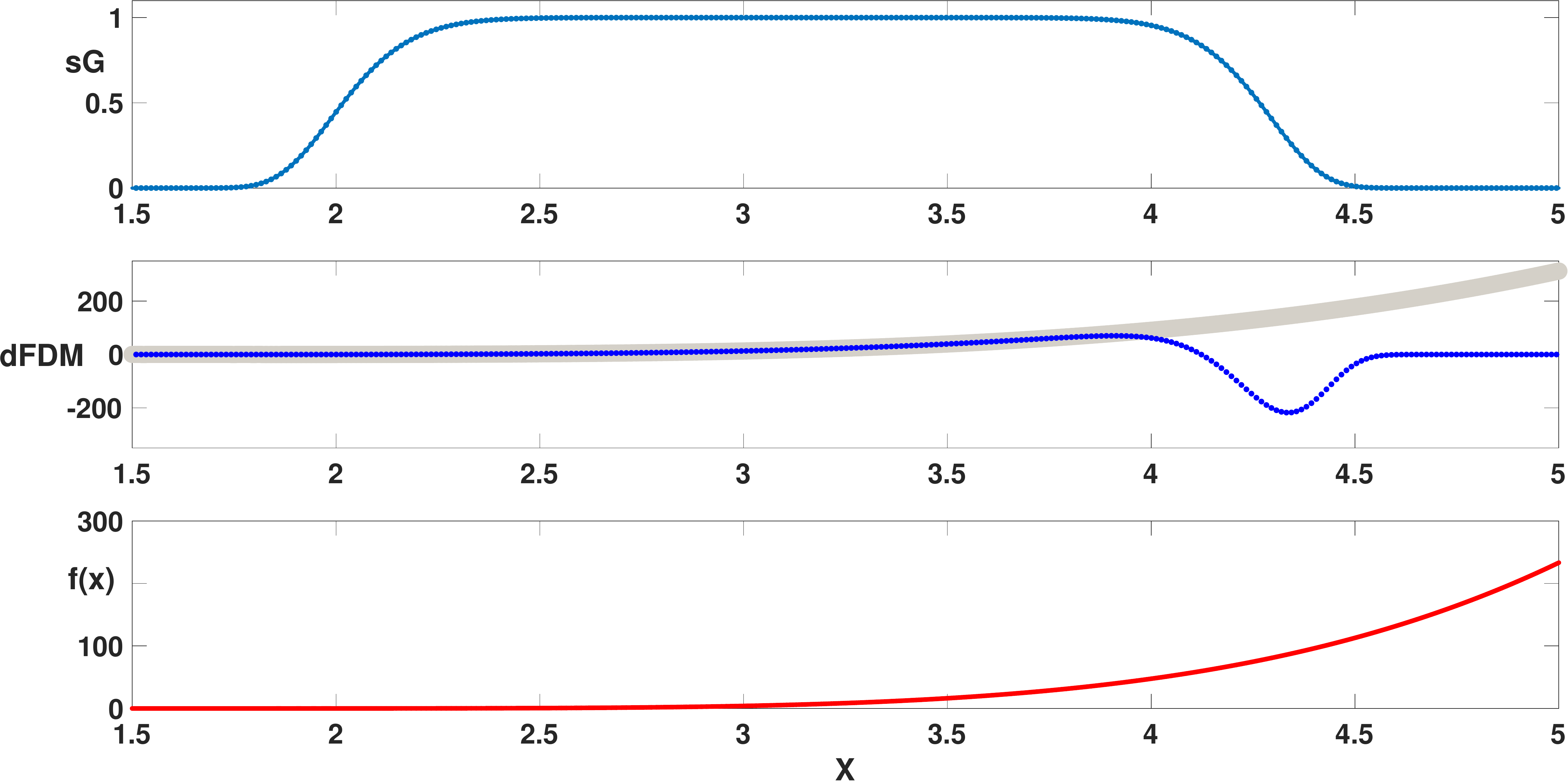}
 \caption{Top: the super-Gaussian $sG(x)$ centered at $x=\pi$, our point of interest, with $s=10$ and $\sigma =1.6$. We have zoomed into the region of interest.
 Middle: the dFDM derivative appears in dots ($N=512$). The exact derivative is displayed in gray. The agreement at $x=\pi$ is up to 14 digits of 
 accuracy. Bottom: the function $f(x)$.}
 \label{sGauss10_512}
 \end{figure}

\begin{figure}[t]
\center
\includegraphics[height=1.5in,width=4in]{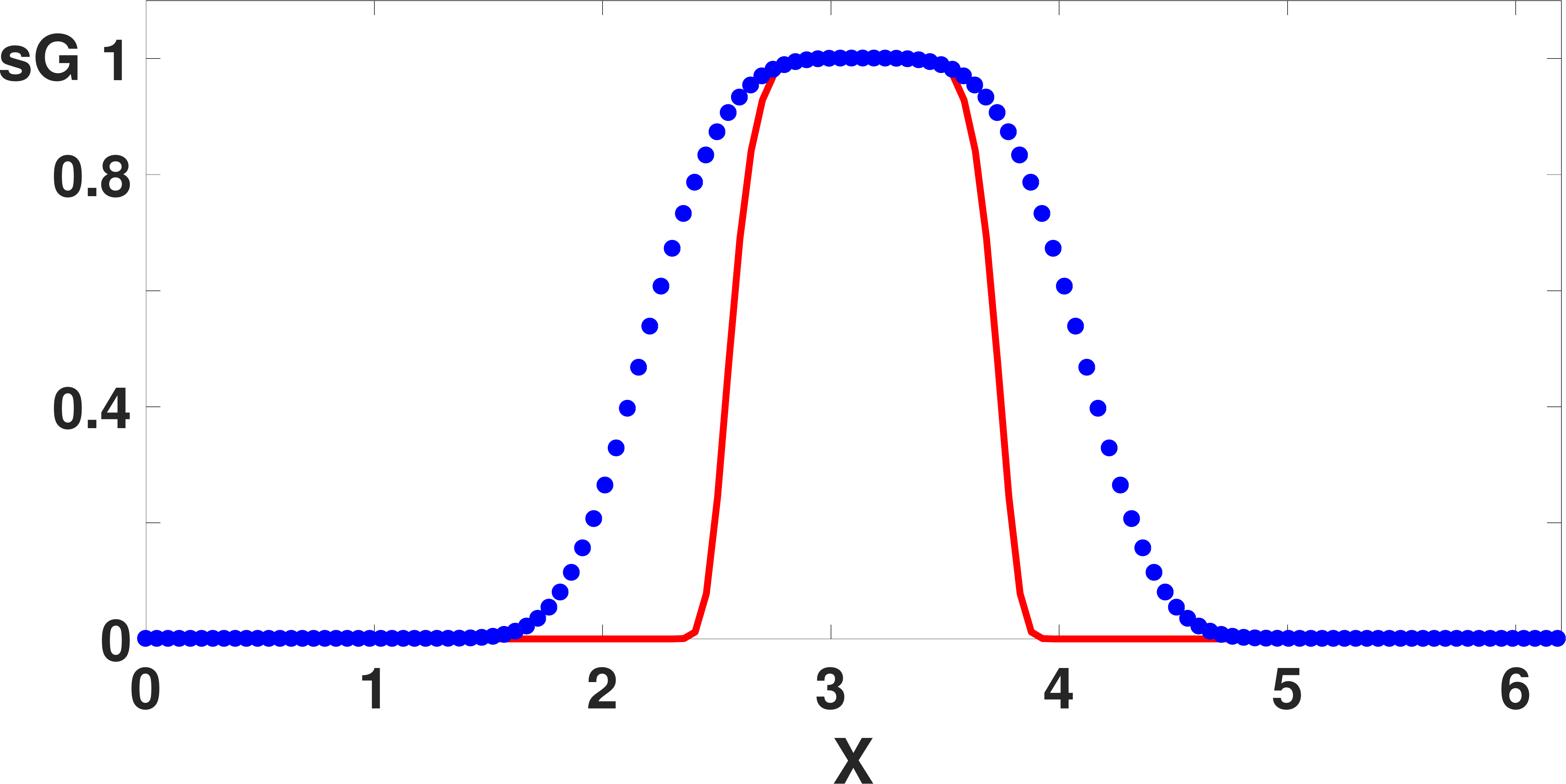}
 \caption{A low resolution example with $N=128$. Two super-Gaussians $sG(x)$ were tested, both centered at $\pi$. 
 The $sG$ depicted with dots used $s=4$ and $\sigma =1.6$ while the $sG$ with a solid line $s=8$ and $\sigma =3.0$.}
 \label{sGauss4_8_128}
 \end{figure}

%
   
 \section{Conclusions}
 ~
 
The main goal of this work is to present a conceptually new framework for a finite difference operator  
  which, as opposed to other finite differences,  has superior spectral accuracy. 
 The method is based solely on numerical manipulations of  singular integrals. 
 The singular integral is a convolution with a distribution. The distribution sets the stage for a multi-resolution finite difference operator.
 The numerical accuracy of the distributional finite difference scheme (dFDM) is illustrated through a series of examples. These are compared
 with the Fourier (FFT based) spectral method. The FFT  method performs differentiation in the frequency domain, whereas the dFDM in
 the physical domain. It is shown that their roundoff error structure is quite different. The dFDM is also compared with a complex-step method
 where the goal is to make the truncation error as small as possible by taking a very small step size $h$. 
 By construction our truncation error is small for smooth functions, which is a property related to  spectral accuracy. We get many digits of accuracy
 with much larger step sizes.
 
 As themes for future investigation, we consider analyzing other singular integral strategies, the possibility of  the making the scheme more 
 compact and efficient 
 computationally while exploring applications to differential equations.
   
 The author's work was supported in part by CNPq under (PQ-1B) 301949/2007-7 and FAPERJ Cientistas do Nosso
Estado project no. 102.917/2011.

\vspace{.5cm}

\noindent{\sc  IMPA, Instituto Nacional de Matem\'atica Pura e Aplicada, Rio de Janeiro, Brasil.}

\noindent{\it E-mail: nachbin@impa.br}

\end{document}